\documentclass{amsart}
\newtheorem{Theorem}{Theorem}[section]
\newtheorem{Definition}[Theorem]{Definition}
\newtheorem{Proposition}[Theorem]{Proposition}

\newtheorem{Lemma}[Theorem]{Lemma}
\newtheorem{Corollary}[Theorem]{Corollary}
\theoremstyle{remark}
\newtheorem{Example}[Theorem]{Example}

\def\eps{\varepsilon}

\def\ovr{\overline}

\def\Gm{\Gamma}

\def\dl{\delta}

\def\bd{\partial}
\def\lm{\lambda}

\def\si{\sigma}

\def\sm{\setminus}
\def\sbs{\subset}

\def\wtl{\widetilde}

\def\be{\begin{enumerate}}
\def\ee{\end{enumerate}}
\def\bT{\begin{Theorem}}
\def\eT{\end{Theorem}}
\def\bP{\begin{Proposition}}
\def\eP{\end{Proposition}}
\def\bD{\begin{Definition}}
\def\eD{\end{Definition}}
\def\bE{\begin{Example}}
\def\eE{\end{Example}}
\def\bL{\begin{Lemma}}
\def\eL{\end{Lemma}}
\def\bC{\begin{Corollary}}
\def\eC{\end{Corollary}}

\def\rS{{\mathcal S}}

\begin{document}
\title{Stein neighborhoods of graphs of holomorphic mappings}
\author{Evgeny A. Poletsky}
\begin{abstract} In this paper we provide sufficient conditions for the graphs
of holomorphic mappings on compact sets in complex manifolds to have
Stein neighborhoods. We show that under these conditions the mappings
have properties analogous to properties of holomorphic functions on
compact sets in $\mathbb C^n$.
\end{abstract}
\thanks{The author was supported by the NSF Grant DMS-0900877.}
\address{Department of Mathematics,  Syracuse University, \newline
215 Carnegie Hall, Syracuse, NY 13244} \email{eapolets@syr.edu}
\maketitle
\section{Introduction}
\par Let $N$ and $M$ be complex manifolds and let $K$ be a compact set in $N$.
There are several different ways to define holomorphic mappings between
$K$ and $M$ but whatever definition is chosen there are at least three
questions we would like to answer:\be\item Can any such mapping be
approximated by holomorphic mappings on neighborhoods of $K$?
\item If two such mappings are close, are they homotopic to each other
    in the space of holomorphic mappings on $K$?
\item Is there a flexible way to shift these mappings?
\ee
\par These questions are either well-studied or trivial when $M$ and $N$ are
Euclidean spaces. One can easily mimic this situation when $K$ and its
image have Stein neighborhoods by imbedding these neighborhoods into
Euclidean spaces. But, in general, neither $K$ nor its image need to
have Stein neighborhoods. For example, they may contain compact complex
submanifolds. A more exquisite example of this kind due to Rosay, when a
mapping of a closed unit ball into an Euclidean space is an embedding
of the open unit ball, can be found in \cite{For}.
\par If we define holomorphic mappings on $K$ as restrictions of
holomorphic mappings on its neighborhoods and assume that $K$ has a
basis of Stein neighborhoods, then we can look at the graph of the
mapping on such a neighborhood $U$. Since the graph is a Stein complex
submanifold of $U\times M$ by Siu's theorem (\cite{S}) we see that the
graph has a Stein neighborhood. After that we can proceed as above.
\par But if we give more reasonable definitions for holomorphic mappings on
compact sets, then the situation is much more subtle. In \cite{DrnFor}
and \cite{For} one can find an excellent discussion of arising problems
and there solutions in many cases. However, finally, everything is
reduced to the question: When the graph of a holomorphic mapping has a
basis of Stein neighborhoods? And the goal of this paper is to provide
a partial but usable answer to this question and also the answers to
three questions above.
\par The author is grateful to Franc Forstneri\v{c} for fruitful and
inspiring discussions.
\section{Fusion}
\par In the future we will need to fuse plurisubharmonic functions.
This means that given a finite collection of open sets $U_j$ and
plurisubharmonic functions $u_j$ on $U_j$ we want to find a
plurisubharmonic function $u$ on the union $U$ of all $U_j$ whose
values are comparable with values of $u_j$. In general, fusion is
impossible. For example, one can cover the infinite hyperplane in
$\mathbb {CP}^n$ by open sets and define plurisubharmonic functions on
these sets equal to $-\infty$ exactly on the hyperplane, but there are
no plurisubharmonic functions defined on a neighborhood of the infinite
hyperplane and equal to $-\infty$ on the hyperplane.
\par Analyzing the proof of Theorem 1 about approximations
of plurisubharmonic functions in \cite{FW}  P. Gauthier in \cite{Gau} came to a way to fuse two (pluri)subharmonic functions on domains in $(\mathbb
C^n)\mathbb R^n$.  This method was modified by N. Gogus in \cite{Go} to
fuse finitely many functions. The main ingredient of the proof in
\cite{Gau} was the existence of a strictly (pluri)subharmonic function
on $(\mathbb C^n)\mathbb R^n$.
\par We use this idea to find a more delicate way of fusing. For this we
introduce a {\it fusing device} which consists of a finite collection
of triples $(V_j,U_j,\chi_j)$, $1\le j\le k$, where $V_j\sbs\sbs U_j$
are open sets in a complex manifold $N$ and $\chi_j$ is a
$C^{\infty}$-function on $N$ with compact support in $U_j$, equal to 1
on $V_j$ and taking values between 0 and 1 elsewhere and  a strictly
plurisubharmonic function $\dl$ on $N$ taking values between 0 and 1.
Let us take the maximal $\lm>0$ such that the functions
$w_j=\dl+\lm\chi_j$ are plurisubharmonic on $N$ for all $j$. We define
the fusing norm as $c=3+2\lm^{-1}$.
\par The advantage of our fusing result is that the same fusing device
can fuse plurisubharmonic functions defined on different neighborhoods.
\bT\label{T:fuse} Let $P$ be a holomorphic mapping of a complex manifold $M$
into a complex manifold $N$. Suppose that $\{\{(V_j,U_j,\chi_j), 1\le
j\le k\},\dl\}$ is a fusing device on $N$ and $c$
is its fusing norm.
\par Let $K$ be a compact set in $M$ such that
$K'=P(K)\sbs\cup_{j=1}^kV_j$. Suppose also that there are open sets
$W_j\sbs N$, $1\le j\le k$, such that $U_j\cap K'\sbs W_j$ and
$K'\sbs\cup_{j=1}^kW_j$ and plurisubharmonic functions $u_j$ on $\wtl
W_j=P^{-1}(W_j)$ never taking $-\infty$ as their values. Then there are
a neighborhood $Z\sbs\cup_{j=1}^k W_j$ of $K'$ and a plurisubharmonic
function $u$ on $\wtl Z=P^{-1}(Z)$ of $K$  such that
\[\|u-u_j\|_{\wtl W_j\cap\wtl Z}\le c\max\{\|u_l-u_m\|_{\wtl W_l\cap
\wtl W_j},1\le m\le l\le k\}.\]
\eT
\begin{proof}  Let
\[b=\max\{\|u_l-u_m\|_{\wtl W_l\cap \wtl W_j},1\le m\le l\le k\}.\]
If $b=0$ then we can take $u(z)=u_j(z)$ if $z\in\wtl W_j$. If $b=\infty$
then we can take as $u$ any plurisubharmonic function on $M$, for example,
$\dl\circ P$. So we assume that $0<b<\infty$.
\par Let $v_j(z)$ be functions on $M$ equal to $u_j(z)+2\lm^{-1}bw_j(P(z))$
on $\wtl W_j\cap\wtl U_j$, where $\wtl U_j=P^{-1}(U_j)$, and to $-\infty$
elsewhere. On $\wtl W=\cup_{j=1}^k(\wtl W_j\cap\wtl U_j)$ we define the
function
\begin{equation}\label{e:du}
u(z)=\max\{v_j(z):\, 1\le j\le k\}.
\end{equation}
\par Let us show that $u$ is plurisubharmonic on a neighborhood $\wtl Z$ of
$K$. If $z\in\wtl W$ does not belong to the boundary of any set $\wtl
W_j\cap\wtl U_j$, then $u$ is plurisubharmonic near $z$.
\par Now suppose that $z\in K$ and the set $L=\{j_1,\dots,j_p\}$ of indexes
such that $z\in\bd(\wtl W_{j_l}\cap\wtl U_{j_l})$ for some $j_l\in L$
is non empty. If $j_l\in L$ then $\zeta=P(z)\in\bd U_{j_l}$ because
$K'\cap U_{j_l}\sbs W_{j_l}$. Hence, $\chi_{j_l}\equiv0$ on the
intersection of some neighborhood $Q_{j_l}$ of $\zeta$ with
$W_{j_l}\cap U{j_l}$.
\par There is $j$ such that $\zeta=P(z)\in V_j\cap K'$. We may assume that
$\chi_j\equiv 1$ on $Q_{j_l}$. Consequently, for
$w\in P^{-1}(Q_{j_l}\cap W_{j_l}\cap U_{j_l})$ and $\xi=P(w)$ we have
\[v_{j_l}(w)-v_j(w)=u_{j_l}(w)-u_j(w)-
2b(\chi_j(\xi)-\chi_{j_l}(\xi))\le-b<0.\] Hence $v_j(w)>v_{j_l}(w)$.
\par Now let $Q'$ be the intersection of all $Q_{j_l}$, $j_l\in L$, and
let $Q\sbs Q'$ be a neighborhood of $\zeta$ such that $Q\cap\bd(W_l\cap
U_l)=\emptyset$ for $l\not\in L$. If $z\in P^{-1}(Q)$ then
\[u(z)=\max\{v_j(z):\,j\not\in L\}\] and this means that
$u$ is plurisubharmonic on $P^{-1}(Q)$. Therefore, there is  a neighborhood
$Z\sbs\cap_{j=1} W_j$ of $K'$ such that $u$ is
plurisubharmonic on $\wtl Z$.
\par If $z\in\wtl W_j\cap\wtl Z$ then there is $l$ between 1 and $k$ such that
\[u(z)-u_j(z)=v_l(z)-u_j(z)=u_l(z)-u_j(z)+2b(\lm^{-1}\dl(P(z))+\chi_j(P(z))).\]
Hence
\[\|u-u_j\|_{\wtl W_j\cap\wtl Z}\le b(3+2\lm^{-1}).\]
\end{proof}
\par This theorem can be used to construct strictly pseudoconvex
functions on $M$. Let us recall (see \cite{Nar}) that an upper
semicontinuous function $u$ on a complex manifold $M$ is called {\it
strictly plurisubharmonic} if for every $z\in M$ there is a
neighborhood $V$ and a $C^\infty$ strictly plurisubharmonic function
$v$ such that $u-v$ is plurisubharmonic on $V$.
\bC\label{C:spsh} Suppose, additionally, that in assumptions of Theorem
\ref{T:fuse} the fibers $P^{-1}(z)$ are manifolds for all $z\in N$ and
the functions $u_j$ are smooth and strictly plurisubharmonic on
$P^{-1}(z)\cap\wtl W_j$ for all $z\in W_j$. Then there is a
neighborhood $V$ of $K$ and a strictly plurisubharmonic function on
$V$.
\eC
\begin{proof} For the proof we note that the functions $v_j$ in the proof of
Theorem \ref{T:fuse} are smooth and and strictly plurisubharmonic on
$P^{-1}(z)\cap\wtl W_j$ for every $z\in\wtl W_j$. Therefore, the
functions $v_j(z)+\dl(P(z))$ are strictly plurisubharmonic on $\wtl
W_j$. If a function $u$ is defined by (\ref{e:du}) then
\[u(z)+\dl(P(z))=\max\{v_j(z)+\dl(P(z)):\, 1\le j\le k\}\]
is strictly plurisubharmonic on a neighborhood $\wtl Z$ of $K$ as the maximum of finitely many strictly plurisubharmonic functions.
\end{proof}
\section{Stein neighborhoods}
\par Let $K$ be a compact set in a complex manifold $N$ and let $M$ be a
complex manifold. We denote by $A(K,M)$ the set of all continuous
mappings of $K$ into $M$ which are holomorphic in the interior of $K$.
Let $H(K,M)$ be the set of all continuous mappings of $K$ into $M$
which can be uniformly approximated with any precision by holomorphic
mappings of neighborhoods of $K$ into $M$ and let $H_{loc}(K,M)$ be the
set of all continuous mappings of $K$ into $M$ with the following
property: for any $f\in H_{loc}(K,M)$ and any point $z\in K$ there is a
neighborhood $U$ of $z$ such that $f$ can be uniformly approximated
with any precision on $K\cap\ovr U$ by mappings into $M$ holomorphic on
neighborhoods of $K\cap\ovr U$. Clearly, $H(K,M)\sbs H_{loc}(K,M)\sbs
A(K,M)$.
\par The famous Mergelyan theorem states that if $N=\mathbb C$ then
$H(K,\mathbb C)=A(K,\mathbb C)$.  In several variables such a theorem
does not hold, for example, when $K$ is a closed disk in $\mathbb C^2$
since $K$ has no interior. Until a better definition for $A(K,M)$ will
show up an analog of the Mergelyan theorem can be searched only when
$K$ is the closure of its interior (see \cite{GoP}).
\par Following the classical theory we continue to say that $A(K,M)$ has the {\it Mergelyan property} if
$A(K,M)=H(K,M)$ and introduce the {\it local Mergelyan
property} claiming that $A(K,M)=H_{loc}(K,M)$. Clearly, $A(K,M)$ has the local
Mergelyan property if $A(K,\mathbb C)$ has it. For example, if $K$ is
the closure of any domain in $N$ with $C^1$-boundary or any compact set
in a Riemann surface, then $A(K,M)$ has the local Mergelyan property.
\par The problem of Stein neighborhoods for graphs of holomorphic
mappings is well discussed in \cite{For}, where it was proved that if
$K$ is the closure of a strongly pseudoconvex domain in $M$ and $f\in
A(K,M)$, then the graph of $f$ on $K$ has a basis of Stein
neighborhoods. A reader can also find there counterexamples which show
that such a statement fails in general settings.
\par In this paper Forstneri\v{c} poses the following problem: Let $K$ be a
compact set with a Stein neighborhood basis in a complex manifold $N$.
Assume that $f\in H(K,M)$.  Does the graph of $f$ on $K$ admits a basis
of open Stein neighborhoods in $N\times M$? The following theorem
gives an affirmative answer to this question even in more general form.
\bT\label{T:snloc} Suppose that a compact set $K$ in a complex manifold $N$
has a basis of Stein neighborhoods. If $M$ is a
complex manifold of dimension $m$ and $f\in H_{loc}(K,M)$, then the
graph of $f$ on $K$ has a basis of Stein neighborhoods in $N\times M$.
\eT
\begin{proof} If $f:\,A\to B$ is a mapping and $C\sbs A$, then $\Gm^f_C$
will denote the graph of $f$ on $C$.
\par The proof will follow in several steps. At the first step we select
an appropriate covering of $\Gm^f_K$. For this we pick up a Riemann
metric $\rho$ on $M$. For every point $z\in M$ there is its neighborhood
$B_z$, a biholomorphic mapping $\phi_z$ of $B_z$ onto the unit ball in
$\mathbb C^m$ and a constant $a_z>1$ such that
\[a_z^{-1}\|\phi_z(z_1)-\phi_z(z_2)\|\le\rho(z_1,z_2)\le
a_z\|\phi_z(z_1)-\phi_z(z_2)\|.\]
\par Then we take a neighborhood $X'\sbs N\times M$ of $\Gm_K^f$ and
find a Stein neighborhood $Y\sbs\sbs N$ of $K$ and a continuous
extension $\wtl f$ of $f$ to $\ovr Y$ such that $\Gm^{\wtl f}_{\ovr Y}\sbs X'$
and there is a strictly plurisubharmonic function $\dl$ on $Y$ with
values between 0 and 1 on $Y$. We cover $\wtl f(\ovr Y)$ by finitely
many neighborhoods $B_{z_j}=B_j$, $1\le j\le k$, and let
$\phi_j=\phi_{z_j}$ and $a=\max_{1\le j\le k} a_{z_j}$. If $z_1,z_2\in
B_j\cap B_l$ then
\[a^{-2}\|\phi_j(z_1)-\phi_j(z_2)\|\le\|\phi_l(z_1)-\phi_l(z_2)\|\le
a^2\|\phi_j(z_1)-\phi_j(z_2)\|.\]
\par For every point $\zeta\in K$ there are its neighborhood $U_\zeta\sbs Y$
and  a set $B_j=B_{j\zeta}$ such that $\Gm^{\wtl f}_{U_\zeta}\sbs\sbs
U_\zeta\times B_j\sbs X'$. Moreover, we may assume that $f$ can be
uniformly approximated with any precision on $K\cap\ovr U_\zeta$ by
mappings holomorphic on neighborhoods of $K\cap\ovr U_\zeta$. We choose
a domain $V_\zeta\sbs\sbs U_\zeta$ containing $\zeta$ and find a finite
covering of $K$ by open sets $V_j=V_{\zeta_j}$, $1\le j\le p$. We set
$U_j=U_{\zeta_j}$ and $B_j=B_{j\zeta}$.
\par Let $U=\cup_{j=1}^pU_j$. Now we choose $r>0$ so small that for every
$1\le j\le p$ the open set \[\wtl U_j=\{(\zeta,z):\,\zeta\in
U_j,\rho(\wtl f(\zeta),z)<r\}\sbs\sbs U_j\times B_j.\] Let
$X=\{(\zeta,z):\,\zeta\in U,\rho(\wtl f(\zeta),z)<r\}$ and let
$P(\zeta,z)=\zeta$ be a projection of $X$ onto $U$. Clearly,
$\Gm_K^f\sbs X\sbs X'$.
\par Our next step is to choose appropriate plurisubharmonic functions
on $\wtl U_j$. For this we fix $\eps>0$ whose precise value will be
determined later. For each $1\le j\le p$ we take a holomorphic mapping
$g_j$ of a neighborhood of $\ovr U_j\cap K$ into $M$ such that
$\rho(g_j(\zeta),f(\zeta))<\eps$ on $\ovr U_j\cap K$. There is a
neighborhood $W_j\sbs Y$ of $\ovr U_j\cap K$ such that $g_j$ is defined
on $W_j$ and $\rho(g_j(\zeta),\wtl f(\zeta))<\eps$ on $W_j$ and $g_j$
maps $W_j$ into $B_j$. Let
$u_j(\zeta,z)=\max\{\log\|\phi_j(z)-\phi_j(g_j(\zeta))\|,\log\eps\}$ be
the functions on $\wtl W_j=P^{-1}(W_j)$.
\par The functions $u_j$ have the following properties. If
$(\zeta,z)\in\bd X$ and $\rho(z,\wtl f(\zeta))=r$, then
\[u_j(\zeta,z)\ge\log(a^{-1}\rho(z,g_j(\zeta)))\ge\log(a^{-1}(r-\eps)).\]
Suppose that $(\zeta,z)\in\wtl U_j\cap\wtl U_l$. Then
\begin{equation}\begin{aligned}
&\|\phi_l(z)-\phi_l(g_l(\zeta))\|\le \|\phi_l(z)-\phi_l(g_j(\zeta))\|+
\|\phi_l(g_j(z))-\phi_l(g_l(\zeta))\|\notag\\&
\le\|\phi_l(z)-\phi_l(g_j(\zeta))\|+2a\eps
\le a^2\|\phi_j(z)-\phi_j(g_j(\zeta))\|+2a\eps.\notag
\end{aligned}\end{equation}
Thus if $\|\phi_j(z)-\phi_j(g_j(\zeta))\|<\eps$, then
\[u_l(\zeta,z)-u_j(\zeta,z)\le\log\eps+\log(a^2+2a)-\log\eps=\log(a^2+2a)=b.\]
If $\|\phi_l(z)-\phi_l(g_l(\zeta))\|<\eps$, then
\[u_l(\zeta,z)-u_j(\zeta,z)\le\log\eps-\log\eps=0.\]
If $\|\phi_j(z)-\phi_j(g_j(\zeta))\|\ge\eps$ and
$\|\phi_l(z)-\phi_l(g_l(\zeta))\|\ge\eps$, then
\[u_l(\zeta,z)-u_j(\zeta,z)\le\log(a^2+2a)=b.\] Thus
\[\|u_j-u_l\|_{\wtl U_j\cap\wtl U_l}\le b.\]
Finally, if $\zeta\in W_j$ then $u_j(\zeta,\wtl
f(\zeta))\le\log(a\eps)$.
\par By Theorem \ref{T:fuse} there is a constant $c>0$ which
does not depend on $\eps$, a neighborhood $Z$, which we may assume to be Stein,
of $K$ and a plurisubharmonic function $u$ on $\wtl Z=P^{-1}(Z)$ such that
\[\|u-u_j\|_{\wtl W_j\cap\wtl Z}\le cb.\]
\par Now if $\zeta\in Z$ and $\rho(z,f(\zeta))=r$, then
\[u(\zeta,z)\ge \log(a^{-1}(r-\eps))-cb,\]
while if $z=\wtl f(\zeta)$ then
\[u(\zeta,z)\le cb+\log(a\eps).\] So if we take $\eps>0$ so small that
$\log(a^{-1}(r-\eps))>\log(a\eps)+2cb$, then the open set
\[D=\{(\zeta,z)\in\wtl Z:\, u(\zeta,z)<\log(a^{-1}(r-\eps))\}\]
contains $\Gm^{\wtl f}_Z$ and $u(\zeta,z)=\log(a^{-1}(r-\eps))$ when
$(\zeta,z)\in\bd D$ and $\zeta\in Z$. Hence there is a plurisubharmonic
function $v$ on $D$ such that $v(\xi,z)\to\infty$ as $(\xi,z)$
approaches a point $(\zeta,z)$, $(\zeta,z)\in\bd D$ and $\zeta\in Z$.
\par Finally, we take a plurisubharmonic exhaustion function $\psi$ on $Z$ and
let $v(\zeta,z)=u(\zeta,z)+\psi(\zeta)$. This is a plurisubharmonic
exhaustion function on $D$.
\par A theorem of R. Narasimhan (see \cite{Nar} and \cite{AnNar}) claims
that a manifold is Stein if it has a plurisubharmonic exhaustion
function and a strictly plurisubharmonic function. Let us construct the
latter function on $X$. For this on each $\wtl W_j$ we take the
functions $u_j(\zeta,z)=\|\phi_j(z)\|^2$ and apply Corollary
\ref{C:spsh}. This implies that $D$ is a Stein manifold.
\par A Stein neighborhood of $K$ and the parameter $r$ can be
chosen as small as we want and this shows that $K$ has a basis of Stein
neighborhoods in $N\times M$.
\end{proof}
\par This theorem has an immediate corollary.
\bC\label{C:sn} If in assumptions of Theorem \ref{T:snloc} $f\in A(K,M)$ and,
additionally, $A(K,M)$ has the local Mergelyan property, then the graph
of $f$ on $K$ has a basis of Stein neighborhoods in $N\times M$.
\eC
\section{Applications}
\par Throughout this sections we assume that $M$ is a complex manifold with a Riemann metric $\rho_M$ and $K$ is a compact set  with a basis of Stein neighborhoods in a complex manifold $N$.
\par The first  application allows us to construct Stein
neighborhoods for special products.
\bT\label{T:spsn} If $K'=K\times[0,1]\sbs N'=N\times\mathbb C$ and
$f(z,t), z\in K, t\in[0,1],$ is a continuous mapping of $K'$ into $M$
such that $f(\cdot,t)\in H(K,M)$ for each $t\in[0,1]$, then $f\in
H_{loc}(K',M)$ and $\Gm^f_{K'}$ has a basis of Stein neighborhoods in
$N'=N\times\mathbb C\times M$.
\eT
\begin{proof} Clearly, $K'$ has a basis of Stein neighborhoods in
$N\times\mathbb C$. By Theorem \ref{T:snloc} all we need to prove is
that $f\in H_{loc}(K',M)$. For this let us take any point $z_0\in K$
and consider the compact set $K'_0=\{z_0\}\times[0,1]\sbs N'$. Clearly,
$H_{loc}(K'_0,M)=C(K'_0,M)$ and, therefore, by Theorem \ref{T:snloc}
the graph $\Gm^f_{K'_0}$ has a Stein neighborhood $Y$ in $N'\times M$.
Let us denote by $Y_1$ a neighborhood of $\Gm^f_{K'_0}$ which compactly
belongs to $Y$. There is a neighborhood $V$ of $z_0$ in $N$ such that
$\Gm^f_{V'}\sbs\sbs Y_1$, where $V'=V\times[0,1]$.
\par Let $F$ be an imbedding of $Y$ into ${\mathbb C}^p$ as a complex
submanifold. By \cite[Theorem 8.C.8]{GR} there are an open neighborhood
$U$ of $F(\ovr Y_1)$ in ${\mathbb C}^p$ and a holomorphic retraction
$P$ of $U$ onto $F(Y)$. We set $\hat f(z,t)=F(z,t,f(z,t))$ for
$(z,t)\in\ovr V'$.
\par Let us show that the function $\hat f$ can be uniformly approximated
with any precision on $\ovr V'\cap K'$ by holomorphic mappings defined
on the neighborhoods of $\ovr V'\cap K'$. For this for $\eps>0$ and for
each $t\in[0,1]$ we can find $\eps_t$ between $0$ and $\eps$ and a
neighborhood $Z_t$ of $K$ and a mapping $f_t\in A(Z_t,M)$ such that
$\rho_M(f_t(z),f(z,t))<\eps_t$ when $z\in K$ and $(z,t,f_t(z))\in Y$
when $z$ is in a neighborhood $W_t$ of $\ovr V\cap K$. For $z\in W_t$ we set
$\hat f_t(z)=F(z,t,f_t(z))$. There is a constant $C$ depending only on
$Y_1$ and $F$ such that $\|\hat f_t(z,t)-\hat f(z,t)\|<C\eps$ when
$z\in W_t\cap K$.
\par Thus we can select points $0=s_0<s_1<\dots<s_q=1$ and $t_j\in[s_{j-1},s_j]$,
$1\le j\le q$ such that $\|\hat f_{t_j}(z)-\hat f(z,t)\|<2C\eps$ when
$t\in[s_{j-1},s_j]$ and $z\in W\cap K$. Set $W=\cap_{j=1}^qW_{t_j}$. We
define a mapping $\hat g$ of $W'=W\times[0,1]$ into $\mathbb C^p$ as
follows: in the first step we set $\hat g(z,t)=\hat f_{t_j}(z)$ when
$t\in[s_{j-1},s_j)$ and in the second step we choose intervals
$[x_j,y_j]$ around point $s_j$, $j=1,\dots,q-1$, and set
\[\hat g(z,t)=\frac{y_j-t}{y_j-x_j}\hat f_{t_{j-1}}(z)+\frac{t-x_j}{y_j-x_j}\hat f_{t_j}(z)\]
on $[x_j,y_j]$. Choosing these intervals small we can require that
$\|\hat g(z,t)-\hat f(z,t)\|<3C\eps$.
\par Since $\hat g$ is a continuous mapping of $W'$ we can use the
Stone--Weierstrass theorem to approximate it uniformly on $\ovr V'\cap
K'$ by holomorphic mappings $g_l$ which are polynomials of degree at most $2l$
in $t$. For example we can use the convolution formula (51) from the
proof of Theorem 7.26 in \cite{Ru}. The obtained functions $\hat g_l$
will be defined on $W\times\mathbb C$ and the direct verification of
the proof in \cite{Ru} shows that they holomorphic in $z$ and converge
uniformly to $\hat g$ as $l\to\infty$.
\par Picking up $\eps>0$ sufficiently small and the approximation above
sufficiently precise we can be assured that $\hat g_l$ maps some neighborhood of $\ovr V'\cap K'$  into
$Y$ and $\|\hat f(z,t)-P(\hat g_l(z,t))\|<4C\eps$ on $\ovr V'\cap K'$.
\par  Let $g_l=\Pi\circ F^{-1}\circ P\circ\hat g_l$, where $\Pi$ is the
projection $N\times M$ onto $M$. Due to the continuity of
$\Pi\circ F^{-1}$ choosing $\eps$ small we guarantees that $g_l$
uniformly approximates $f$ with any prescribed precision.
\end{proof}
\par In the results below we will use polynomials taking prescribed
values at given points. While on $\mathbb C$ this job is done by
Lagrange interpolating polynomials, there is no a canonical procedure
in several variables. To avoid ambiguity we will use the following
simple result.
\bL\label{L:ip} Suppose that a finite set
$\si\sbs\mathbb C^q$ of distinct points $\{\si_1,\dots,\si_q\}$ and a
set $h=\{h_1,\dots,h_q\}\sbs\mathbb C$. Let $\|h\|=\max_{1\le j\le
q}|h_j|$. Then for every $R>0$ there is a constant $C$ depending only
on $\si$ such that there is a polynomial $L_{\si,h}$ of degree at most
$q$ on $\mathbb C^q$ taking values $h_j$ at all points $\si_j$ and
whose uniform norm does nor exceed $C\|h\|$ when $\|z\|\le R$.
\eL
\begin{proof} Let us fix $q$ hyperplanes defined as solutions of
the equations $L_j(z)=\langle a_j,z\rangle+b_j=0$ with $\|a_j\|=1$ such
that $L_j(\si_i)=0$ if and only if $i=j$. Set
\[L_{\si,h}(z)=\sum_{k=1}^qh_k\frac{\prod_{j\ne k}L_j(z)}
{\prod_{j\ne k}L_j(\si_k)}.\] It is easy to see that $L_{\si,h}$ has
all required properties.
\end{proof}
\par The following application is the Mergelyan-type
approximation of mappings in $A(K,M)$. When $K\sbs\mathbb C$ this
problem was studied in \cite{Cha}.
\bT\label{T:at} Suppose that $\si=(\si_1,\dots,\si_q\}$ is a finite set in
$K$, $h=\{h_1,\dots,h_q\}$ is a set in $M$ and $A(K,\mathbb C)$ has the
Mergelyan property. For every $f\in H_{loc}(K,M)$ and for every
$\eps>0$ there is $\dl>0$ such that if $\rho_M(f(\si_j),h_j)<\dl$,
$1\le j\le q$, then one can find a neighborhood $U$ of $K$ and $g\in
A(U,M)$ such that $\rho_M(f(\zeta),g(\zeta))<\eps$ for all $\zeta\in K$
and $g(\si_j)=h_j$, $1\le j\le q$.
\eT
\begin{proof} Since $f\in H_{loc}(K,M)$, by Theorem \ref{T:snloc} the
graph $\Gm^f_K$ has a Stein neighborhood $Y$ in $N\times M$. Let $F$ be
an imbedding of $Y$ into ${\mathbb C}^p$ as a complex submanifold. By
\cite[Theorem 8.C.8]{GR} there are an open neighborhood
$U\sbs\sbs\mathbb C^p$ of $F(\Gm^f_K)$ in ${\mathbb C}^p$ and a
holomorphic retraction $P$ of $U$ onto $F(Y)$.
\par Let $\hat f(\zeta)=F(\zeta,f(\zeta))$. Since $A(K,\mathbb C)$ has the
Mergelyan property for any $\eta>0$ there is a neighborhood $V$ of $K$
and $\hat g_1\in A(V,\mathbb C^p)$ such that $\|\hat f(\zeta)-\hat
g_1(\zeta)\|<\eta$ on $K$ and $\hat g_1(V)\sbs\sbs U$. Choose $\dl>0$
such that $\|F(\si_j,f(\si_j))-F(\si_j,h_j)\|<\eta$, $1\le j\le q$.
Then $\|\hat g_1(\si_j)-F(\si_j,h_j)\|<2\eta$. By Lemma \ref{L:ip} the
uniform norm on $U$ of an interpolating polynomial mapping $L:\mathbb
C^p\to\mathbb C^p$ such that $L(\hat g_1(\si_j))=F(\si_j,h_j)-\hat
g_1(\si_j)$  does not exceed $c\eta$ on $U$,
where $c$ depend only on $U$ and $\si$. If $\hat g(\zeta)=\hat
g_1(\zeta)+L(\hat g_1(\zeta))$ then $\hat g(\si_j)=F(\si_j,h_j)$ and $\|\hat
f(\zeta)-\hat g(\zeta)\|<(1+c)\eta$ on $K$.
\par Define $g=\Pi\circ F^{-1}\circ P\circ\hat g$. Due to the continuity of
$\Pi\circ F^{-1}\circ P$ choosing $\eta$ sufficiently small guarantees
that $g$ has all needed properties.
\end{proof}
\par This theorem has an interesting corollary which was proved by E. L.
Stout in \cite{St} when $A(K,\mathbb C)=C(K,\mathbb C)$.
\bC\label{C:ep} If $A(K,\mathbb C)$ has the Mergelyan property then
$A(K,M)$ has the Mergelyan property.
\eC
\begin{proof} If $A(K,\mathbb C)$ has the Mergelyan property then
$A(K,M)$ has the local Mergelyan property and by the previous result
$H_{loc}(K,M)=H(K,M)$. Hence $A(K,M)=H(K,M)$.
\end{proof}
\par To give another application we choose a Riemann metrics $\rho_N$ on
$N$ and introduce a topology on the set $\rS(N,M)=\rS$ of all pairs
$(K,f)$, where $K$ is a compact set in $N$ and $f\in A(K,M)$. Given any
continuous extension $\Phi$ of $f$ to $N$ and $\eps>0$ we define a
$\Phi,\eps$-neighborhood of $(K,f)$ as a set of all pairs $(L,g)\in\rS$
such that the Hausdorff distance $\rho(K,L)$ between $K$ and $L$ is
less than $\eps$  and $\rho_M(\Phi(\zeta), g(\zeta))<\eps$ for all
$\zeta\in L$.
\par It is easy to verify that if $U$ is a $(\Phi,\eps)$-neighborhood
of $(K,f)$, $V$ is a $(\Psi,\dl)$-neighborhood of $(L,g)$ and $(N,h)\in
U\cap V$, then there is a neighborhood of $(N,h)$ lying in $U\cap V$.
Hence our choice of neighborhoods defines a topology on $\rS$.
\par Suppose that we have a family of compact sets $L_t\sbs N$, $0\le t\le1$,
and a finite set $\si\sbs L_t$ for all $t$. Such a family is {\it
Rad\'o continuous with respect to $\si$} if there are homeomorphisms
$\phi_t\in A(L_0,N)$ mapping $L_0$ onto $L_t$ such that for every
$t\in[0,1]$ the set $\si$ belongs to the set of fixed points of
$\phi_t$ and  the mappings $\phi_s\circ\phi^{-1}_t$ converge uniformly
to identity as $s\to t$.
\par In \cite{Ra} T. Rad\'o gave necessary and sufficient conditions
for a family $L_t$ of closures of connected and simply connected
domains in $\mathbb C$ to be Rad\'o continuous. For us this notion is
useful because, evidently, if a family $L_t$ is Rad\'o continuous and
$g_0\in A(L_0,M)$, then the path $(L_t,g_0\circ\phi^{-1}_t)$ is
continuous in the topology introduced above.
\par The following theorem shows that under reasonable assumptions
$\Phi,\eps$-neighborhoods  are ``convex".
\bT\label{T:ht} Suppose that $f\in H(K,M)$. For every $\eps>0$ there
is $\dl>0$ with the following properties: if:\be
\item $\si=\{\si_1,\dots,\si_q\}$ is a set of distinct points in $N$;
\item  $L_t$, $0\le t \le1$, is a family of compact
    sets $L_t\sbs N$ lying in the $\rho_N,\dl$-neighborhood of $K$ such that $\si\sbs L_t\cap K$ when $0\le t\le 1$;
\item the family $L_t$ is Rad\'o continuous with respect to $\si$;
\item mappings $g_0\in A(L_0,M)$ and $g_1\in A(L_1,M)$ and  there is a
    continuous extension $\Phi$ of $f$ to $N$ such that
    $\rho_M(g_0(\zeta),\Phi(\zeta))<\dl$ and
    $\rho_M(g_0(\zeta),\Phi(\zeta))<\dl$ when $\zeta\in L_0$ or
    $\zeta\in L_1$ respectively;
\item $g_0(s_j)=g_1(s_j)$, $1\le j\le q$,
\ee then there is a continuous path $(L_t,g_t)$, $0\le t\le 1$, such that
$\rho_M(g_t(\zeta),\Phi(\zeta))<\eps$ and $g_t(\si_j)=g_0(\si_j)$,
$1\le j\le q$ and $0\le t\le 1$.
\eT
\begin{proof} By Theorem \ref{T:snloc} the graph $\Gm^f_K$ has a Stein
neighborhood $Y$ in $N\times M$. Let $F$ be an imbedding of $Y$ into
${\mathbb C}^p$ as a complex submanifold. By \cite[Theorem 8.C.8]{GR}
there are an open neighborhood $U\sbs\sbs\mathbb C^p$ of $F(\Gm^f_K)$
in ${\mathbb C}^p$ and a holomorphic retraction $P$ of $U$ onto $F(Y)$.
If $\Psi$ is a mapping of $A\sbs N$ into $M$ then we set $\hat
\Psi(\zeta)=F(\zeta,\Psi(\zeta))$, $\zeta\in A$.
\par By Theorem \ref{T:at} for any $\eta>0$ there is a neighborhood
$D$ of $K$ and $h\in A(D,M)$ such that $\|\hat
h(\zeta)-\hat\Phi(\zeta)\|<\eta$ when $\zeta\in D$, $\Gm^h_D\sbs Y$ and
$F(\Gm^{\hat h}_D)\sbs U$. Moreover, we can find $h$ such that
$h(\si_j)=g_0(\si_j)$ provided $\rho_M(f(\si_j),g_0(\si_j))$ doess not
exceed some $\dl>0$ for all $j$.
\par We may assume that $\dl$ is so small that all $L_t\sbs D$. Let $\hat h_0$
and $\hat h_1$ be the restrictions of $\hat h$ to $L _0$ and $L_1$
respectively. Let $\hat d_0(\zeta)=\hat h_0(\zeta)-\hat g_0(\zeta)$ and
$\hat d_1(\zeta)=\hat h_1(\zeta)-\hat g_1(\zeta)$. Then
\[\|\hat d_0(\zeta)\|\le\|\hat h(\zeta)-\hat f(\zeta)\|+
\|\hat f(\zeta)-\hat g_0(\zeta)\|\le\eta+\|\hat f(\zeta)-\hat
g_0(\zeta)\|\] for $\zeta\in L_0$ and a similar inequality holds for
$\hat d_1$.  Reducing if needed the size of $\dl$ we may assume that
$\|\hat d_0\|$ on $L_0$ and $\|\hat d_1\|$ on $L_1$ do not exceed
$2\eta$.
\par On $L_t$ we define
\[\hat b_t(\zeta)=(1-t)\hat d_0(\phi_t^{-1}(\zeta))+t\hat
d_1(\phi_1(\phi^{-1}_t(\zeta))),\] where $\phi_t\in A(L_0,N)$ are
homeomorphisms of $L_0$ onto $L_t$ keeping $\si$ fixed. Clearly,
$\|b_t(\zeta\|<2\eta$ for all $t$ and $\zeta\in L_t$, $\hat b_0=\hat
d_0$, $\hat b_1=\hat d_1$ and $\hat b_t(\si_j)=0$. By an observation
preceding the statement of the theorem the path $(L_t,\hat b_t)$ is
continuous in $\rS(N,\mathbb C^p)$.
\par Let $\hat g_t(\zeta)=\hat h(\zeta)-\hat b_t(\zeta)$. Then $\hat
g_t(\si_j)=\hat g_0(\si_j)$ and \[\|\hat g_t(\zeta)-\hat
\Phi(\zeta)\|\le\|\hat h(\zeta)-\hat \Phi(\zeta)\|+2\eta\le3\eta.\]
\par Define $g_t=\Pi\circ F^{-1}\circ P\circ\hat g_t$. Due to the continuity
of $\Pi\circ F^{-1}\circ P$ choosing $\eta$ sufficiently small
guarantees that $g_t$ has all needed properties.
\end{proof}
\par One of the essential problem in complex geometry is moving objects
preserving their complex structures. The problem was investigated in
\cite{Gr}, where the notion of sprays was introduced. Later this notion
obtained a significant development in papers of Forstneri\v{c} and his
coauthors. We will use the definition of sprays close to one given in
\cite{DrnFor}.
\par Let $f\in A(K,M)$. A spray with an exceptional set $\si=\si(f)\sbs K$
is a map $F : K\times V\to M$, where $V$ (the parameter set of $F$) is
an open subset of an Euclidean space $\mathbb C^N$ containing the
origin, such that the following holds:\be\item $F$ is in $A(K\times
V,M)$;\item $F(z,0)=f(z)$ and $F(z,t)=f(z)$ for $z\in\si$ and $t\in V$;
\item for every $z\in K\sm\si$ the map $\bd_tF(z,t):\mathbb C^N\to
    T_{F(z,t)}M$ is surjective (the domination condition).\ee
\par The following theorem establishing an existence of sprays can be
proved exactly as Lemma 4.2 in \cite{DrnFor} where it was established
for compact sets with smooth boundary. For the sake of completeness we
give another proof here.
\bT\label{T:sp} Let $f\in H_{loc}(K,M)$ and $\si=\{\si_1,\dots,\si_k\}$ be
a finite set in $K$. Then there is a spray $G$ with the exceptional set
$\si$ and $G(z,0)=f(z)$ on $K$.
\eT
\begin{proof} We take a Stein neighborhood $Y$ of $\Gm^f_K$ in $N\times M$
and imbed it into $\mathbb C^p$ as a complex submanifold by a
holomorphic mapping $G$. Let $\hat\si_j=F(\si,f(\si_j))$ and
$\hat\si=\{\hat\si_1,\dots,\hat\si_k\}$. There are exactly $p$
polynomials $P_1,\dots,P_N$ on $\mathbb C^p$ such that
$\hat\si=\{P_1=\dots=P_p=0\}$. This can be proved by induction in $p$.
If $p=1$ then $P_1$ is a polynomial of degree $k$ with simple roots in
$\hat\si$. If it is proved for $p-1$ then we represent $\mathbb C^p$ as
$\mathbb C^{p-1}\times\mathbb C=\{(z',\zeta):\,z'\in\mathbb
C^{p-1},\zeta\in\mathbb C\}$. We may assume that no two points in
$\hat\si$ have the same coordinates in $\mathbb C^{p-1}$. Let
$\hat\si'=\{\hat\si'_1,\dots,\hat\si'_k\}$ be the projection of
$\hat\si$ onto $\mathbb C^{p-1}$. By the induction of assumptions there
are polynomials $P_1,\dots,P_{p-1}$ on $\mathbb C^{p-1}$ such that
$\hat\si'=\{P_1=\dots=P_{p-1}=0\}$. Denote by $\hat\si''_j$ the
projection of $\hat\si_j$ on $\mathbb C$ and take an interpolating
polynomial $Q_p$ on $\mathbb C^{p-1}$ such that
$Q_p(\hat\si'_j)=\hat\si''_j$. Let $P_p(z',\zeta)=\zeta-Q_p(z')$.
Clearly, $\hat\si=\{P_1=\dots=P_p=0\}$.
\par There are an open neighborhood $U$ of $F(\Gm^f_K)$ in ${\mathbb C}^p$ and
a holomorphic retraction $P$ of $U$ on $G(Y)$. We may assume that $U$
is so small that $dP$ has the maximal rank on $U$. There is $\eps>0$ so
small that if all $v_j\in\mathbb C^N$, $1\le j\le N$, have the norm
less than $\eps$ and $v=(v_1,\dots,v_N)$, then
\[H(\zeta,v)=F(\zeta,f(\zeta))+P_1(\zeta,f(\zeta))v_1+\dots+
P_N(\zeta,f(\zeta))v_N\in U\] for all $\zeta\in K$. Let $\Pi$ be the
natural projection of $M\times N$ onto $M$. Hence the mapping
$G(\zeta,v)=\Pi\circ F^{-1}\circ P\circ H(\zeta,v)$ is well defined
when $v$ is in the ball of radius $\eps$ in $\mathbb C^{p^2}$,
holomorphic, $G(\zeta,0)=f(\zeta)$ if $\zeta\in\si$ and for every
$\zeta\in K\sm\si$ the map $\bd_vG(\zeta,v):\mathbb C^{p^2}\to
T_{G(\zeta,v)}M$ is surjective. To show the last statement we assume
that $P_1(\zeta,f(\zeta))\ne0$ and let $v_2=\dots=v_N=0$. Then it is
easy to see that $\bd_{v_1}G(\zeta,v_1):\mathbb C^{p^2}\to
T_{G(\zeta,v_1)}M$ is surjective because $\bd_{v_1}H$, $dP$ and
$d(\Pi\circ(F^{-1})$ are surjective.
\end{proof}

\end{document}